\documentclass{article}

\usepackage{float}
\usepackage{subcaption}
\usepackage{PRIMEarxiv}
\usepackage{amsmath}
\usepackage[utf8]{inputenc} 
\usepackage[T1]{fontenc}    
\usepackage{hyperref}       
\usepackage{url}            
\usepackage{booktabs}       
\usepackage{amsfonts}       
\usepackage{nicefrac}       
\usepackage{microtype}      
\usepackage{lipsum}
\usepackage{fancyhdr}       
\usepackage{graphicx}       
\graphicspath{{media/}}     

\title{A Boundary Integral-based Neural Operator for Mesh Deformation 
\thanks{
\textbf{the code will be available upon request}} 
}

\author{
  Zhengyu Wu\\
  School of Electronics and Information Engineering \\
  Hangzhou Dianzi University \\
  Hangzhou, China\\
  \texttt{23070225@hdu.edu.cn} \\
   \And
  Jun Liu\\
  School of Electronics and Information Engineering \\
  Hangzhou Dianzi University \\
  Hangzhou, China\\
  \texttt{ljun77@hdu.edu.cn} \\
  \And   
  Wei Wang\\
  School of Electronics and Information Engineering \\
  Hangzhou Dianzi University \\
  Hangzhou, China\\
  \texttt{wwang@hdu.edu.cn} \\
}

\begin{document}
\maketitle

\begin{abstract}
This paper presents an efficient mesh deformation method based on boundary integration and neural operators, formulating the problem as a linear elasticity boundary value problem (BVP). To overcome the high computational cost of traditional finite element methods and the limitations of existing neural operators in handling Dirichlet boundary conditions for vector fields, we introduce a direct boundary integral representation using a Dirichlet-type Green's tensor. This formulation expresses the internal displacement field solely as a function of boundary displacements, eliminating the need to solve for unknown tractions. Building on this, we design a Boundary-Integral-based Neural Operator (BINO) that learns the geometry- and material-aware Green's traction kernel. A key technical advantage of our framework is the mathematical decoupling of the physical integration process from the geometric representation via geometric descriptors. While this study primarily demonstrates robust generalization across diverse boundary conditions, the architecture inherently possesses potential for cross-geometry adaptation. Numerical experiments, including large deformations of flexible beams and rigid-body motions of NACA airfoils, confirm the model's high accuracy and strict adherence to the principles of linearity and superposition. The results demonstrate that the proposed framework ensures mesh quality and computational efficiency, providing a reliable new paradigm for parametric mesh generation and shape optimization in engineering.
\end{abstract}


\section{Introduction}

Mesh deformation is a fundamental component of multiphysics simulation, widely applied in electromagnetic reduced-order modeling \cite{Feng2019,Wang2011,Liu2023}, integrated circuit (IC) optimization \cite{baranowski2023design,guo2024optimization}, and complex fluid dynamics scenarios such as wing icing \cite{Morelli2021}, parachute mechanics \cite{takizawa2012computational}, and free-surface flows \cite{fenton1997mesh}. Its primary advantage lies in enabling rapid global mesh updates through boundary node displacements alone, thereby avoiding the computationally expensive process of remeshing while preserving topological consistency \cite{Wang2012}. Existing approaches include linear elasticity \cite{Dwight2009}, Radial Basis Functions (RBF) \cite{DeBoer2007}, Inverse Distance Weighting (IDW) \cite{Witteveen2009}, Laplace smoothing \cite{sorgente2023survey}, FEMWARP \cite{shontz2010analysis}, and spectral deformation \cite{Rong2008}. Among these, the linear elasticity method is widely used. It models deformation as a Dirichlet boundary value problem (BVP) for the linear elastic equation without body forces, typically solved using the Finite Element Method (FEM) \cite{Lamecki2016}. Recently, Physics-Informed Neural Networks (PINNs) \cite{cuomo2022scientific} have emerged as an alternative to solve such partial differential equations (PDEs) and has been used for mesh deformation \cite{Aygun2023,Liu2025}. However, because the training process of a PINN is equivalent to solving the PDE, the network must be retrained whenever boundary conditions or geometries change. This high latency makes PINNs unsuitable for real-time simulation . Therefore, there is a critical need for surrogate models that can generalize across diverse boundary conditions.

To enhance generalization, various neural operators have been developed to learn mappings between function spaces. Notable variants—including Graph Neural Operators (GNO) \cite{Anandkumar2020}, Fourier Neural Operators (FNO) \cite{Li2020}, Wavelet Neural Operators (WNO) \cite{tripura2022wavelet}, and Convolutional Neural Operators (CNO) \cite{raonic2023convolutional}—have improved training speed, convergence, and stability by refining their mathematical architectures \cite{Kovachki2023}. However, traditional neural operator theory primarily addresses parametric partial differential equations (PDEs) with fixed boundaries and source terms. Recently, solving boundary value problems (BVPs) has emerged as a prominent research focus. Some approaches \cite{Kashi} employ a "Lift Operator" to extend boundary values across the entire domain before applying FNO; nevertheless, these methods remain limited to uniform grids and specific boundary regularities. In contrast, Graph Neural Networks (GNNs) can handle irregular domains but often lack physical interpretability \cite{Loetzsch}. Additionally, operator networks such as DeepONet \cite{Lu2019} and its variants \cite{He2024}, including POD-DeepONet \cite{Lu2022}, face limitations: their "Branch Net" architecture typically requires fixed input dimensions and predefined sensor patterns. They lack permutation invariance, meaning their performance deteriorates significantly if the ordering or density of boundary points changes. This limitation reduces their flexibility in handling dynamic mesh deformation tasks involving varying node counts or evolving topologies.

To address these challenges, we propose shifting the computational paradigm from "global fitting" to "boundary-driven" inference. Mathematically, the solution to a linear elastic mesh deformation problem is uniquely determined by its boundary conditions. Recent studies, such as BINet \cite{Lin2021} (Boundary Integral Networks) and BI-GreenNet \cite{Lin2023}, have demonstrated that learning the Boundary Integral Equation (BIE)—rather than the raw partial differential equation (PDE)—enables the reduction of a $d$-dimensional spatial problem to a $(d-1)$-dimensional boundary problem through potential theory. However, these approaches primarily focus on geometric generalization under fixed boundary conditions, rather than variable boundary conditions.

In this paper, we propose the Boundary-Integral-based Neural Operator (BINO), a highly efficient mesh deformation framework that rigorously formulates the task as a linear elastic boundary value problem (BVP). To overcome the computational overhead of the finite element method (FEM) and the limitations of existing neural operators in handling vector-valued Dirichlet conditions, we introduce a direct boundary integral representation based on the Dirichlet-type Green’s tensor. This formulation is advantageous because it expresses the internal displacement field entirely as a function of boundary displacements, thereby mathematically eliminating the need to solve for unknown surface tractions and streamlining the physical modeling process. Building on this foundation, the BINO architecture is designed to learn a Green’s traction kernel capable of capturing both geometric features and material properties. A key technical contribution is the use of geometric descriptors, which achieve a mathematical decoupling of the physical integration process from the geometric representation. This decoupling provides exceptional flexibility: while this study highlights robustness across different boundary conditions, the architecture inherently possesses the potential for cross-topological adaptation, freeing the model from the rigid sensor patterns required by traditional operator networks.

The subsequent sections of this paper are organized as follows. First, the fundamental principles of mesh deformation based on linear elasticity theory are presented. Next, a rigorous derivation of the boundary integral formula is provided using the Dirichlet-type Green's tensor. Within this theoretical framework, we propose a Boundary Integral-based Neural Operator (BINO) and define mesh quality metrics for performance evaluation. The following experimental section details the training data generation process, which includes parametric affine, translational, and harmonic boundary conditions, and outlines the data synthesis and preprocessing methodology. The effectiveness of BINO is validated through numerical case studies involving flexible beams and airfoils, followed by dedicated linear validation. Finally, the research contributions are summarized, and directions for future work are discussed.

\section{Linear Elasticity-based Mesh Deformation}
The mesh deformation problem can be viewed as a form of the linear elasticity problem from structural mechanics. Given displaced boundary nodes, the displacement of remaining nodes can be calculated using linear elasticity equations. The equations can be defined as:
\begin{equation}
    \begin{aligned}
        L^{T} c L u = 0, && u \in D \\
        u = u_{b}, && u \in \partial D
    \end{aligned}
    \label{eq:Linear_elastic}
\end{equation}
where $L$ is the differential operator matrix, $c$ is the material property matrix, $D$ is the interior domain, $\partial D$ is its boundary, $u$ is the displacement vector field to be solved for, and $u_{b}$ is the prescribed displacement function on the boundary. This system represents a homogeneous partial differential equation with Dirichlet boundary conditions. In the two-dimensional case (under the plane strain assumption), The strain-displacement operator $ L$ is defined as:
\begin{equation}
L = 
\begin{bmatrix}
\frac{\partial}{\partial x} & 0 \\
0 & \frac{\partial}{\partial y} \\
\frac{\partial}{\partial y} & \frac{\partial}{\partial x}
\end{bmatrix}
\end{equation}
The constitutive matrix $c$, defined by the elastic modulus $E$ and Poisson’s ratio $\nu$, is given by:
\begin{equation}
c = \frac{E}{(1+\nu)(1-2\nu)}
\begin{bmatrix}
1-\nu & \nu & 0 \\
\nu & 1-\nu & 0 \\
0 & 0 & \frac{1-2\nu}{2}
\end{bmatrix}
\end{equation}
To achieve better mesh quality after deformation, the value of Poisson's ratio $\nu$ should be set between 0.3 and 0.45. Since a high ratio causes element distortion and a low ratio reduces deformation resistance \cite{Shamanskiy2021}.

\section{Boundary Integral Formulation via Dirichlet-type Green’s Tensor}
For a linear elastic, isotropic, and homogeneous body occupying a two-dimensional domain $D \subset \mathbb{R}^2$ with a piecewise smooth boundary $\partial D$, Somigliana's identity provides a boundary integral representation of the displacement field. Under the assumption of plane strain (or plane stress with modified constants), the displacement $u$ at an interior source point $\mathbf{x} \in D$ is expressed as \cite{Brebbia1994}:

\begin{equation}
\label{eq:Somigliana}
u(\mathbf{x}) = \int_{\partial D} U^*(\mathbf{x}, \mathbf{y}) t(\mathbf{y}) \, d\Gamma_{\mathbf{y}} - \int_{\partial D} T^*(\mathbf{x}, \mathbf{y}) u(\mathbf{y}) \, d\Gamma_{\mathbf{y}}
\end{equation}

where $\mathbf{x}$ denotes the source point where a unit concentrated force is applied, $\mathbf{y}$ represents the field point on the boundary $\partial D$, and $d\Gamma_{\mathbf{y}}$ is the boundary line element. The boundary traction vector is defined by $t(\mathbf{y}) = \sigma(\mathbf{y})n(\mathbf{y})$, with $n$ as the unit outward normal at $\mathbf{y}$. 

The displacement field in Eq. \eqref{eq:Somigliana} can be more explicitly expressed in terms of its Cartesian components. Let $U_{il}^*(\mathbf{x}, \mathbf{y})$ and $T_{il}^*(\mathbf{x}, \mathbf{y})$ denote the components of the displacement and traction fundamental solution tensors, respectively. Here, the first subscript $i$ refers to the direction of the response at the field point $\mathbf{y}$, while the second subscript $l$ corresponds to the direction of the unit point force applied at the source point $\mathbf{x}$. For the two-dimensional case, these kernels are defined as \cite{Brebbia1994}:

\begin{equation}
\label{eq:basic_solution_u}
U_{il}^*(\mathbf{x}, \mathbf{y}) = \frac{1}{8\pi G(1-\nu)} \left[ (3-4\nu) \ln \left( \frac{1}{r} \right) \delta_{il} + r_{,i}r_{,l} \right]
\end{equation}

\begin{equation}
\label{eq:basic_solution_t}
T_{il}^*(\mathbf{x}, \mathbf{y}) = -\frac{1}{4\pi(1-\nu)r} \left\{ \frac{\partial r}{\partial n} \left[ (1-2\nu)\delta_{il} + 2r_{,i}r_{,l} \right] + (1-2\nu)(r_{,i}n_l - r_{,l}n_i) \right\}
\end{equation}

where indices $i, l \in \{1, 2\}$, $r = \|\mathbf{y} - \mathbf{x}\|$ is the Euclidean distance, $r_{,i} = (y_i - x_i)/r$ are the direction cosines, and $\partial r / \partial n = r_{,k}n_k$ is the normal derivative. The material constants are the shear modulus $G$ and Poisson's ratio $\nu$. 

In order to calculate the inner displacement field, traditional methods such as the Boundary Element Method (BEM) typically form  boundary integral equations by utilizing the Kelvin fundamental solutions as weighting functions. However, this requires simultaneous knowledge of both boundary tractions $t_i$ and displacements $u_i$. 
In this paper, we use a more efficient approach by employing a Dirichlet-type Green's tensor $G(\mathbf{x}, \mathbf{y})$ specifically tailored for the domain $D$. This tensor is defined so that it satisfies the homogeneous Dirichlet boundary condition \cite{Mazya2013}:
\begin{equation} 
\label{eq:dirichlet_cond}
G(\mathbf{x}, \mathbf{y}) = 0, \quad \forall \mathbf{y} \in \partial D
\end{equation}
By substituting the Kelvin solution $U^*$ in Eq. (\ref{eq:Somigliana}) with the Dirichlet-type Green's tensor $G$ defined in Eq. (\ref{eq:dirichlet_cond}), the first integral term—which involves the boundary tractions $t_i(\mathbf{y})$—vanishes identically. Consequently, the internal displacement field $u_l(\mathbf{x})$ can be determined solely from the prescribed boundary displacements $u_i(\mathbf{y})$ by the following simplified integral representation:
\begin{equation} 
\label{eq:direct_integral}
u(\mathbf{x}) = -\int_{\partial D} T^G(\mathbf{x}, \mathbf{y}) u(\mathbf{y}) d\Gamma_{\mathbf{y}}
\end{equation}
where $T^G(\mathbf{x}, \mathbf{y})$ represents the traction kernel derived from the displacement tensor $G$.This formulation, as shown in Eq. (\ref{eq:direct_integral}), bypasses the conventional requirement of solving a system of linear equations for unknown boundary tractions, effectively transforming the boundary value problem into a direct integration procedure. 

\section{Boundary-Integral-based Neural Operators}
The implementation of the direct integral representation in Eq. (\ref{eq:direct_integral}) relies on the fact that the traction kernel $T^G$ is uniquely determined by the Dirichlet-type Green's displacement tensor $G$. Specifically, $T^G$ represents the boundary traction field resulting from the stress state associated with $G$. Once the explicit form of $G$ is established to satisfy the boundary constraint in Eq. (\ref{eq:dirichlet_cond}), the corresponding traction kernel is obtained by applying the differential operators of elasticity and projecting the resulting stress tensor onto the boundary normal. This one-to-one correspondence ensures that the kernel $T^G$ is a fundamental property of the domain’s geometry and material constants. 

Therefore, we can employ a neural network to approximate the kernel $T^G$ as a geometry and material-aware parameterized function. For a domain characterized by a geometry descriptor $\mathcal{D}$ (which can be instantiated as latent shape embeddings, boundary point clouds, or signed distance fields) and material property $\mathcal{A}$, the neural operator $\mathcal{N}_{\theta}$ learns a generalized mapping from the coordinate pair $(\mathbf{x}, \mathbf{y})$ to the components of the traction tensor:
\begin{equation}
\label{eq:neural_kernel}
T^G(\mathbf{x}, \mathbf{y}) \approx \mathcal{N}_{\theta}(\mathbf{x}, \mathbf{y}; \mathcal{D}; \mathcal{A})
\end{equation}
By conditioning the network on $\mathcal{D}$ and $\mathcal{A}$, the model effectively encapsulates the fundamental mapping between a domain’s manifold and its corresponding Green’s traction kernel. In the present implementation, while $\mathcal{D}$ and $\mathcal{A}$ is kept constant to rigorously evaluate the network’s generalization across diverse boundary conditions $u_b$, the architectural design remains intrinsically compatible with geometric and material variations. By substituting this neural-approximated kernel into Eq. (\ref{eq:direct_integral}), the displacement field at any internal point $\mathbf{x} \in D$ is evaluated through a neural boundary integral:
\begin{equation}
\label{eq:neural_integration}
u(\mathbf{x}) \approx -\int_{\partial D} \mathcal{N}_{\theta}(\mathbf{x}, \mathbf{y}; \mathcal{D};\mathcal{A}) u_b(\mathbf{y}) d\Gamma{\mathbf{y}}
\end{equation}
In this framework, $\mathcal{N}_{\theta}$ serves as a domain-agnostic neural kernel provider. In contrast to conventional Neural Operators that learn mappings between fixed-grid function spaces, the proposed formulation decouples the physical integration process from the geometric and material representation. This separation enables the model to exploit the dimensionality reduction inherent in the Boundary Element Method (BEM). Consequently, the architecture can adapt to diverse geometries and materials by updating the descriptor $\mathcal{D}$, the integration manifold $\partial D$ and the material property $\mathcal{A}$ , while preserving the integrity of the underlying neural logic.

For practical implementation in mesh deformation tasks, the continuous boundary integral in Eq. (\ref{eq:neural_integration}) is discretized using the pre-defined boundary nodes of the initial mesh. Let the boundary $\partial D$ be represented by a set of $K$ discrete points $\{\mathbf{y}^{(k)}\}_{k=1}^K$. The neural boundary integral is thus approximated by the following discrete summation:
\begin{equation}
\label{eq:discrete_sum}
u(\mathbf{x}) \approx - \sum_{k=1}^{K} \mathcal{N}_{\theta}(\mathbf{x}, \mathbf{y}^{(k)}; \mathcal{D}; \mathcal{A}) u_b(\mathbf{y}^{(k)}) \Delta W^{(k)}
\end{equation}
where $\mathbf{y}^{(k)}$ and $u_b(\mathbf{y}^{(k)})$ are the coordinates and prescribed displacements of the $k$-th boundary node, respectively. The integration weight $\Delta W^{(k)}$ is defined as:
\begin{equation}
\Delta W^{(k)} = \frac{L^{(k)}}{C}
\end{equation}
where $L^{(k)}$ denotes the length of the boundary segment associated with node $\mathbf{y}^{(k)}$. The scaling constant $C$ is introduced to normalize the kernel's contribution and improve the conditioning of the neural approximation. Based on our empirical numerical experiments, a value of $C=50$ is recommended.

From a computational perspective, this formulation offers significant efficiency advantages. Unlike traditional volume-based neural operators that scale with the total number of domain nodes, the evaluation of Eq. (\ref{eq:discrete_sum}) for any interior point $\mathbf{x}$ entails a complexity of $O(K)$, where $K$ is the number of boundary nodes. For a mesh containing $M$ interior points, the total complexity of a single inference pass is $O(M \times K)$. Given that $K \ll M$ in most practical engineering geometries, this boundary-integral-based approach effectively reduces the computational overhead compared to fully-connected graph neural networks or dense volumetric operators.

\section{Mesh quality metric}
To quantitatively evaluate the geometric integrity of the triangular mesh after deformation, we employ a shape quality metric $Q$ based on the equilateral relative factor. This metric is defined as \cite{Knupp1999}:
\begin{equation}
Q = 4\sqrt{3} \frac{A}{\sum_{i=1}^{3} l_i^2}
\end{equation}
where $A$ denotes the area of the triangular element and $l_i$ represents the length of its $i$-th edge. Historically rooted in finite element analysis, this formulation characterizes the deviation of a given triangle from an ideal equilateral shape. The normalization constant $4\sqrt{3}$ ensures that $Q$ is bounded within the range $[0, 1]$, where $Q=1$ corresponds to a perfect equilateral triangle and $Q \to 0$ indicates a degenerate or "sliver" element with vanishing area. By utilizing this metric, we can rigorously assess the robustness of the neural boundary integral across different boundary conditions and verify that the predicted displacement fields do not lead to unphysical mesh distortions.

\section{Experiment}
This section evaluates the accuracy and robustness of the proposed mesh deformation framework through three distinct numerical investigations. First, a flexible beam within a rectangular domain is analyzed to assess performance under structural distortion. Subsequently, the framework is applied to a NACA 0012 airfoil involving combined translation and moderate angle rotation to examine its stability under complex rigid-body motions. Furthermore, the linearity and superposition characteristics of the framework are rigorously verified by applying varied displacement scaling factors and multi-load combinations. This test ensures that the surrogate model—a Deep Residual Network (ResNet) featuring four residual blocks with 156 units each—accurately preserves the mathematical structure of the underlying linear elasticity equations. The model is trained using the Adam optimizer ($1 \times 10^{-3}$ initial learning rate) with a Mean Squared Error (MSE) loss function to achieve high-fidelity predictions across both geometric and physical transformations.

\subsection{Training Data Generation}

To ensure the surrogate model effectively learns the linear elastic response under diverse loading scenarios, a comprehensive dataset is constructed using a high-fidelity numerical solver. The data generation process involves prescribing stochastic boundary conditions on the computational domain $\Omega$ and solving the governing partial differential equations (PDEs) via the Finite Element Method (FEM). Two distinct deformation categories—affine-translational transformations and harmonic distortions—are designed to cover both global structural movements and local geometric perturbations.

\subsubsection{Parametrization of Boundary Conditions}
The boundary of the domain $\partial \Omega$ is partitioned into a fixed boundary $\Gamma_{\text{fixed}}$ and a loading boundary $\Gamma_{\text{load}}$. For each sample, the displacement field $\mathbf{g}$ on $\Gamma_{\text{load}}$ is synthesized through the following two mechanisms:

\paragraph{Affine and Translational Transformations} 
To capture fundamental structural responses such as stretching, shearing, and rigid-body translation, a stochastic affine map is applied. For a boundary point $\mathbf{x} = [x, y]^T$, the prescribed displacement $\mathbf{g}(\mathbf{x})$ is defined as:
\begin{equation}
    \mathbf{g}(\mathbf{x}) = \mathbf{M} (\mathbf{x} - \mathbf{x}_c) + \mathbf{b}
\end{equation}
where $\mathbf{x}_c$ denotes the geometric center of the loading boundary. The translation vector $\mathbf{b}$ is sampled from a uniform distribution within the target magnitude range. The affine matrix $\mathbf{M} \in \mathbb{R}^{2 \times 2}$ is generated by:
\begin{equation}
    \mathbf{M} = \frac{\mathbf{M}_{rand}}{\max(|\mathbf{M}_{rand}|)} \cdot A_{target}
\end{equation}
where $\mathbf{M}_{rand}$ is a random matrix with elements in $[-0.5, 0.5]$, and $A_{target}$ controls the deformation scale to remain within the linear elastic regime.

\paragraph{Harmonic Distortion Patterns}
To enhance the model's robustness against complex stress concentrations and localized irregularities, a polar-coordinate-based harmonic perturbation is introduced:
\begin{equation}
    \mathbf{g}(\theta) = \begin{bmatrix} \alpha A \sin(k_1 \theta + \phi_1) \\ \alpha A \cos(k_2 \theta + \phi_2) \end{bmatrix}
\end{equation}
where $\theta = \operatorname{atan2}(y-y_c, x-x_c)$ is the polar angle relative to the center, $k \in [k_{min}, k_{max}]$ represents the spatial frequency, $\phi$ is the random phase, and $\alpha$ is a scaling coefficient (typically $0.1$). This approach ensures the dataset encompasses a wide spectrum of boundary curvatures.

\subsubsection{Data Synthesis and Preprocessing}
For each realization, the prescribed displacement functions are implemented as Dirichlet boundary conditions. The displacement field $\mathbf{u}$ over the entire domain is then computed by solving the equilibrium equations of linear elasticity. To maintain data integrity, a verification step is performed to prune samples that exhibit non-convergence or mesh entanglement.

The final dataset $\mathcal{D} = \{(\mathcal{G}_i, \mathcal{U}_i)\}_{i=1}^{N}$ consists of the boundary responses $\mathcal{G}$ (input) and the corresponding full-field displacements $\mathcal{U}$ (ground truth). Additionally, auxiliary geometric features, including nodal coordinates $\mathbf{X}$, normal vectors $\mathbf{N}$, and integration weights $w$, are stored to facilitate the computation of physics-informed loss functions during the training phase. The dataset is randomly shuffled. We set the Poisson's ratio to 0.3 and the Young's modulus to 2000 to generate high-fidelity FEM training data.

In this study, the training configurations for the two case studies are tailored to their respective verification objectives. For the flexible beam, the training set comprises 1,000 samples, consisting of 500 realizations of random affine transformations and 500 realizations of harmonic distortion patterns, which are shuffled to ensure the model captures both global trends and localized irregularities. For the NACA 0012 airfoil, to specifically strengthen the model's ability to represent rigid-body and affine modes, the training set consists of 1,000 pure affine transformation samples.

Notably, the magnitude of displacements in the training sets is maintained at a standard scale, which is approximately one order of magnitude smaller than that used in the subsequent testing examples. Despite this discrepancy, the model demonstrates high prediction accuracy across all cases. This capability, further corroborated by the linearity and scale-invariance tests in Section \ref{subsec:linearity}, confirms that the surrogate model does not merely interpolate within the training range but successfully learns the underlying mapping patterns. Consequently, the proposed method can generalize to boundary conditions of arbitrary magnitudes while remaining strictly consistent with the principles of linear elasticity.

\subsection{Flexible Beam }
\begin{figure}[H]
    \centering
    \begin{subfigure}[t]{0.48\textwidth}
        \centering
        \includegraphics[width=\textwidth]{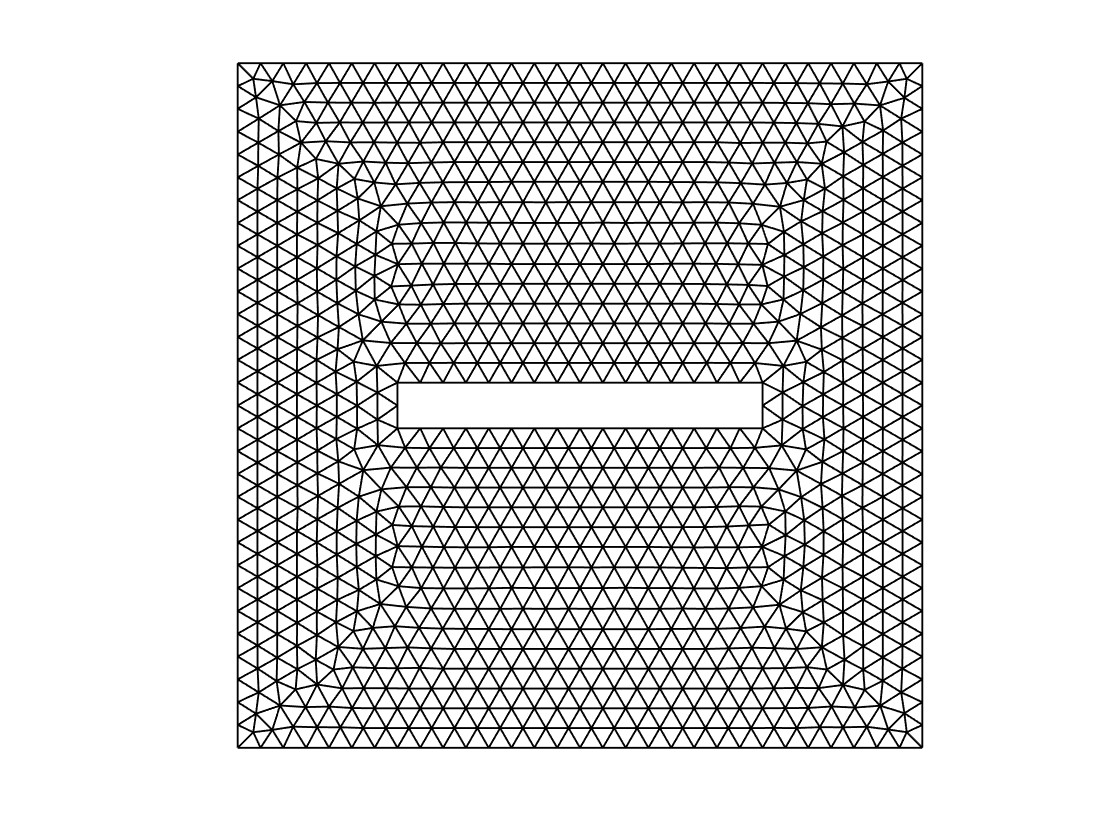}
        \caption{}
        \label{fig:BeamOriginalMesh}
    \end{subfigure}
    \hfill
    \begin{subfigure}[t]{0.48\textwidth}
        \centering
        \includegraphics[width=\textwidth]{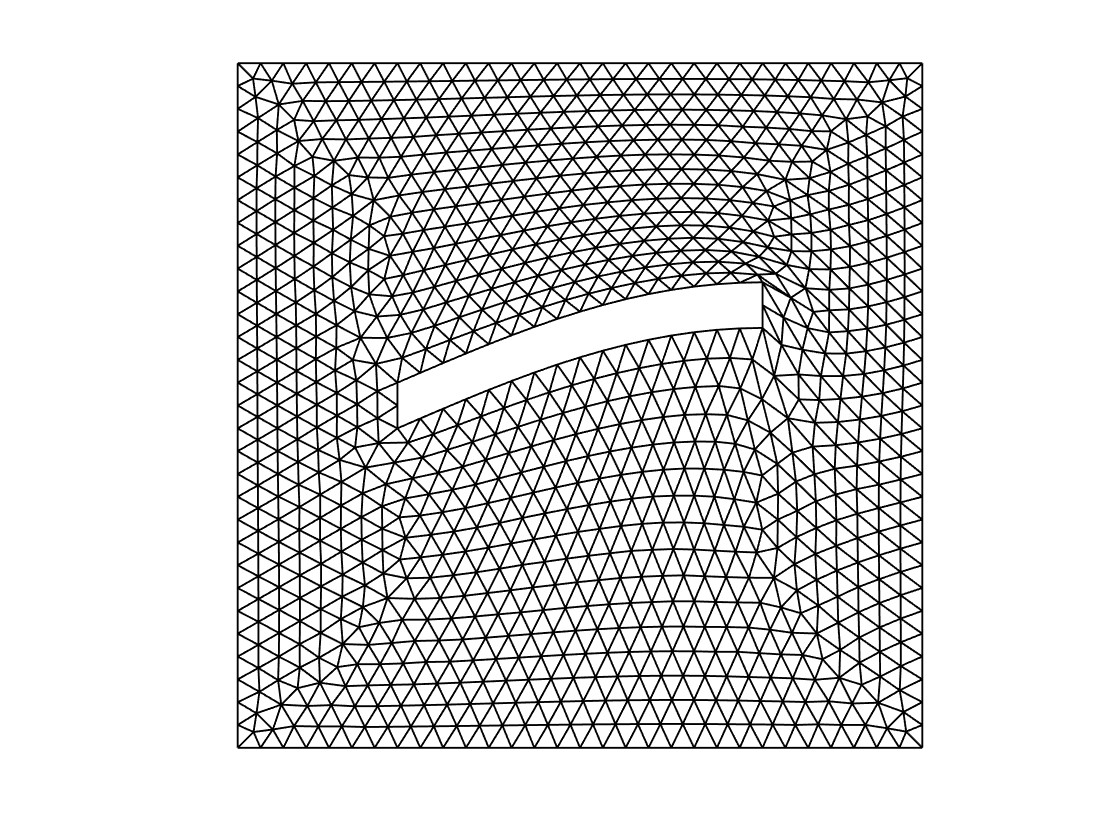}
        \caption{}
        \label{fig:BeamDeformedlMesh}
    \end{subfigure}
    
    \caption{Deformation of Flexible beam : (a) initial mesh (1980 triangles).   (b) deformed mesh (mesh topology and the number of
triangles remain the same).}
    \label{fig:BeamMesh}
\end{figure}

The first case study focuses on a rectangular flexible beam \cite{Shamanskiy2021}.
The initial computational mesh, consisting of 1,980 triangular elements, is shown in Fig. \ref{fig:BeamOriginalMesh}. The beam is fixed on its left end and sits in the center of the domain. The dimensions of the domain are (0, 15)× (0, 15) and the position of the structure is (3.5, 11.5) × (7, 8). The deformation is based on a sinusoidal
function sin($\frac{\pi(x-3.5)}{16}$) with varying amplitude.
Upon applying displacement boundary conditions, the beam exhibits moderate bending, as illustrated in Fig. \ref{fig:BeamDeformedlMesh}.

\begin{figure}[H]
    \centering
    \begin{subfigure}[t]{0.48\textwidth}
        \centering
        \includegraphics[width=\textwidth]{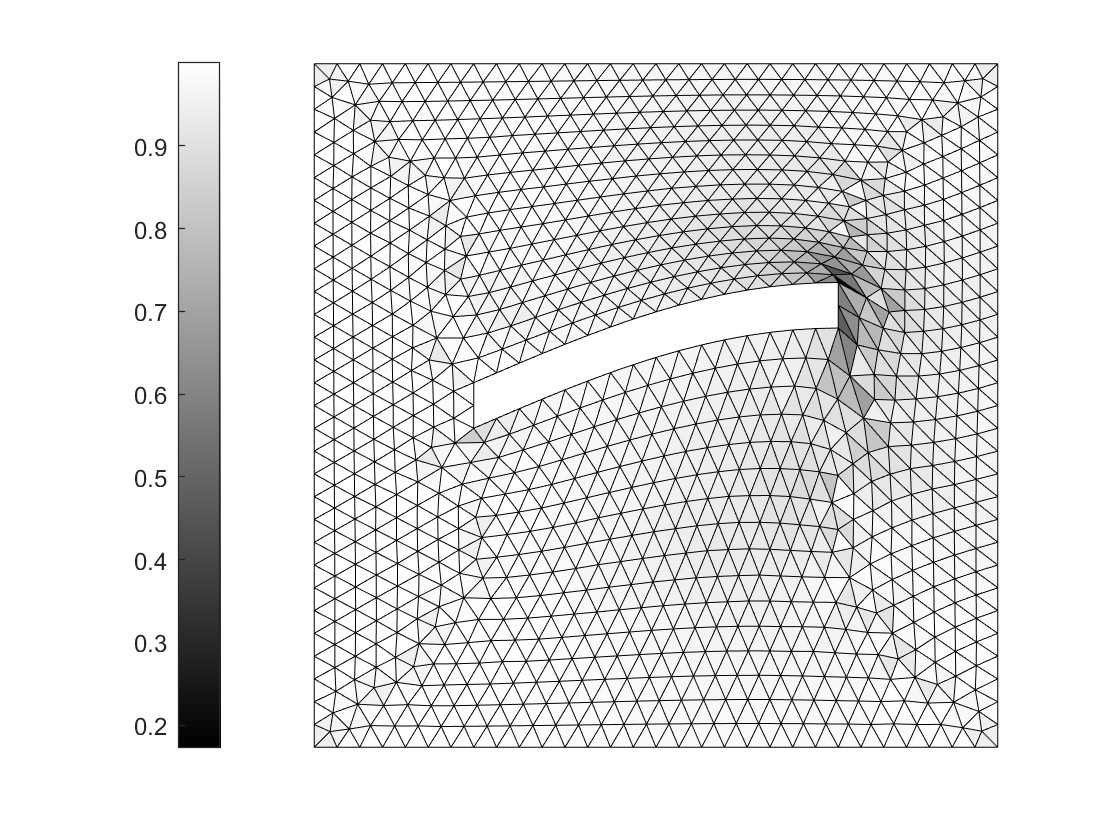}
                \caption{}
        \label{fig:BeamMeshQualityTrue}
    \end{subfigure}
    \hfill
    \begin{subfigure}[t]{0.48\textwidth}
        \centering
        \includegraphics[width=\textwidth]{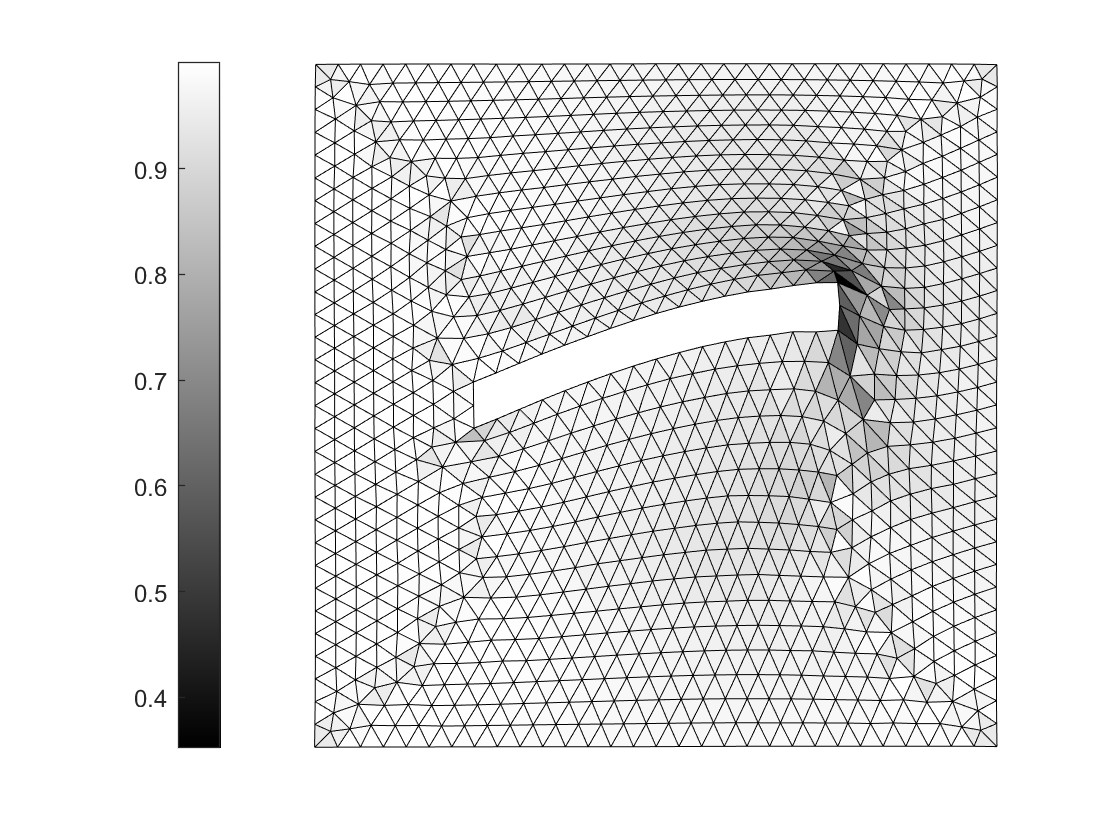}
        \caption{}
        \label{fig:BeamMeshQualityPredict}
    \end{subfigure}

    \caption{Element quality metrics of flexible beam distorted from the top and bottom layers of beam: (a) the deformed mesh generated by FEM. (b) the deformed mesh generated by the proposed surrogate-model.}
\end{figure}

To further assess the geometric fidelity, Figure \ref{fig:BeamMeshQualityTrue} and \ref{fig:BeamMeshQualityPredict} compare the element quality metrics between the traditional Finite Element Method (FEM) and the proposed surrogate model. The results indicate that the surrogate model effectively captures the mesh distortion at the top and bottom layers without introducing degenerate elements.

\begin{figure}[H]
    \begin{subfigure}[t]{1\textwidth} 
        \centering
        \includegraphics[width=1\textwidth]{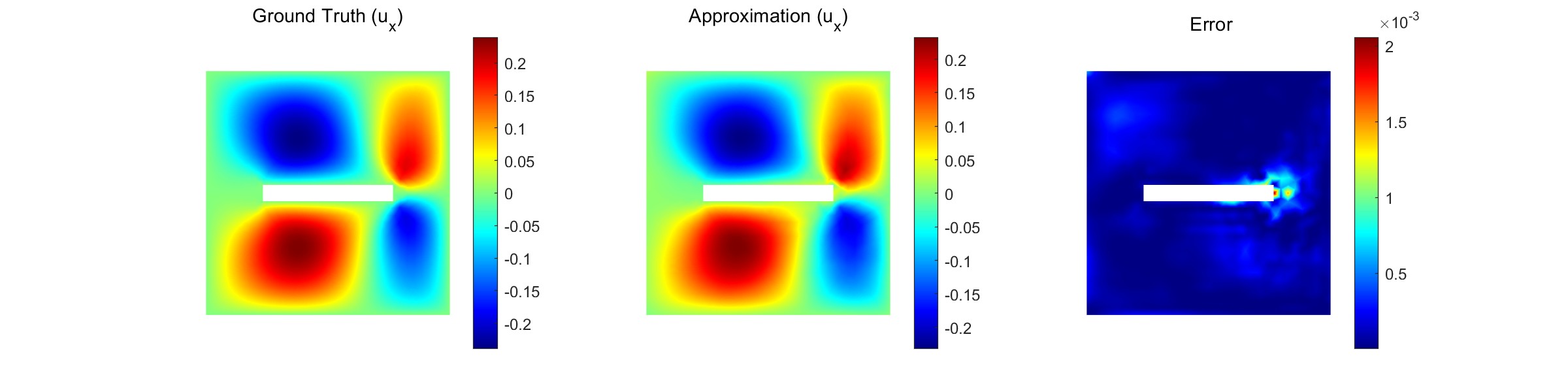} 
        \caption{}
        \label{fig:BeamUx}
    \end{subfigure}

    \begin{subfigure}[t]{1\textwidth}
        \centering
        \includegraphics[width=1\textwidth]{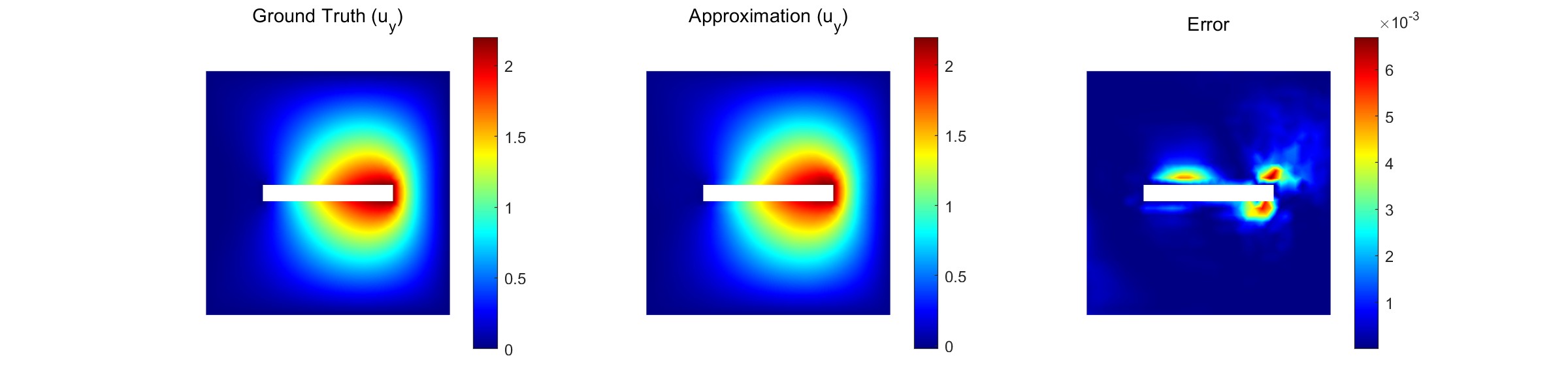} 
        \caption{}
        \label{fig:BeamUy}
    \end{subfigure}
    \caption{The proposed surrogate model for the solution of Eq. (\ref{eq:Linear_elastic}) in this case : (a) the solution in $X$ direction. (b) the solution in $Y$ direction. Both evaluate the approximation accuracy by the point-wise squared error.}
\end{figure}

Finally, Figures \ref{fig:BeamUx} and \ref{fig:BeamUy} illustrate the predicted displacement fields in the $x$ and $y$ directions, respectively, with the approximation accuracy represented by the point-wise squared error. The proposed surrogate model achieves the \textbf{global relative error of only 2.99\%}, demonstrating high-fidelity performance in solving the linear elasticity problem defined in Eq. (\ref{eq:Linear_elastic}). This result highlights the model's potential for real-time structural analysis with accuracy comparable to that of the finite element method (FEM).
However, it should be noted that local variations were observed at the two corner points on the left side of the beam. These minor inaccuracies are primarily attributed to inherent numerical singularities in the boundary integral formulas at points of rapid geometric change, a common phenomenon in boundary-based numerical methods. Despite these local effects, the overall predictive accuracy throughout the entire domain remains robust.

\subsection{Airfoil}

\begin{figure}[H]
    \centering
    \begin{subfigure}[t]{0.48\textwidth}
        \centering
        \includegraphics[width=\textwidth]{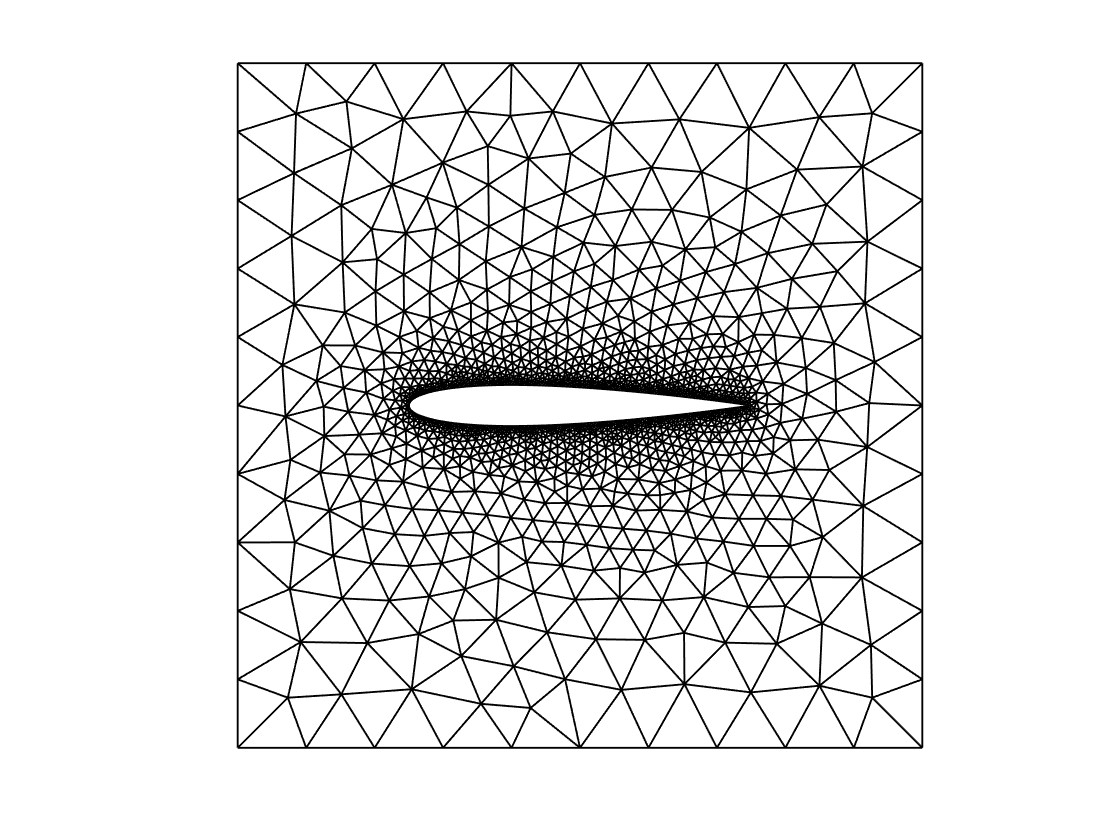}
        \caption{}
        \label{fig:NCAA0012OriginalMesh}
    \end{subfigure}
    \hfill
    \begin{subfigure}[t]{0.48\textwidth}
        \centering
        \includegraphics[width=\textwidth]{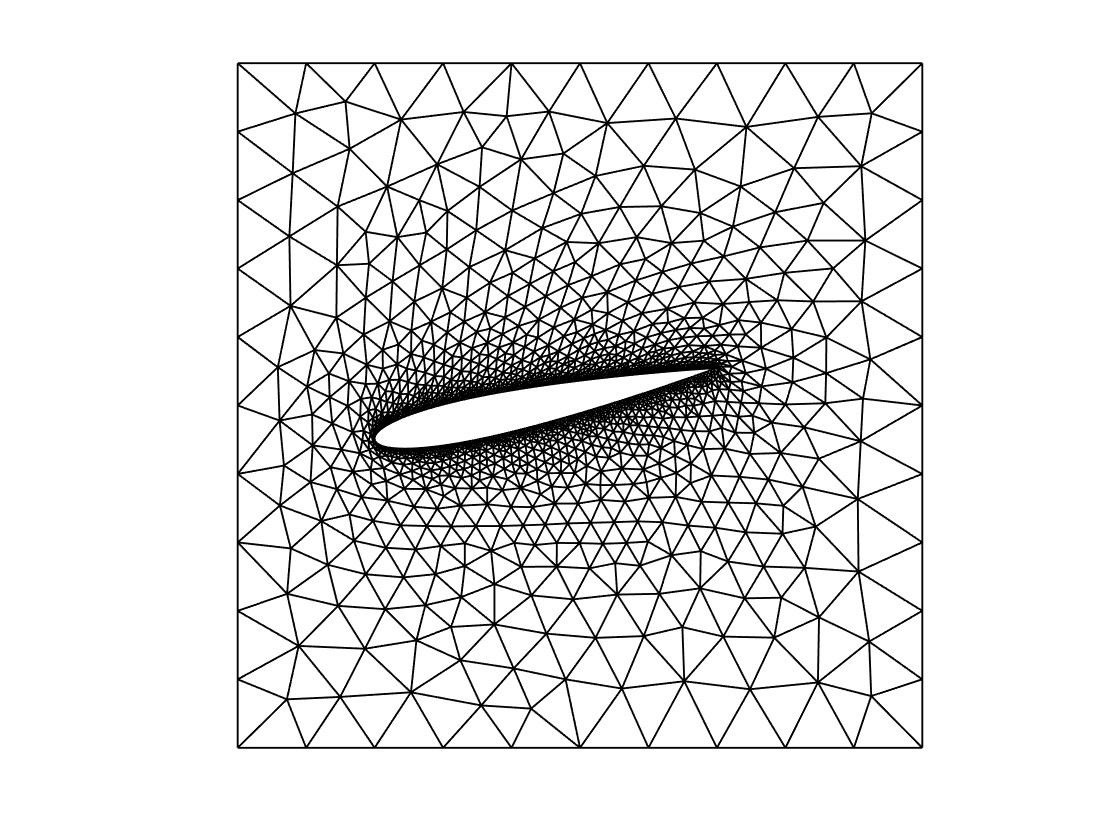}
        \caption{}
        \label{fig:NCAA0012DeformedMesh}
    \end{subfigure}
    \caption{Deformation of NACA 0012 airfoil : (a) initial mesh (2603 triangles). (b) deformed mesh.}
\end{figure}

\begin{figure}[H]
    \centering
    \begin{subfigure}[t]{0.48\textwidth}
        \centering
        \includegraphics[width=\textwidth]{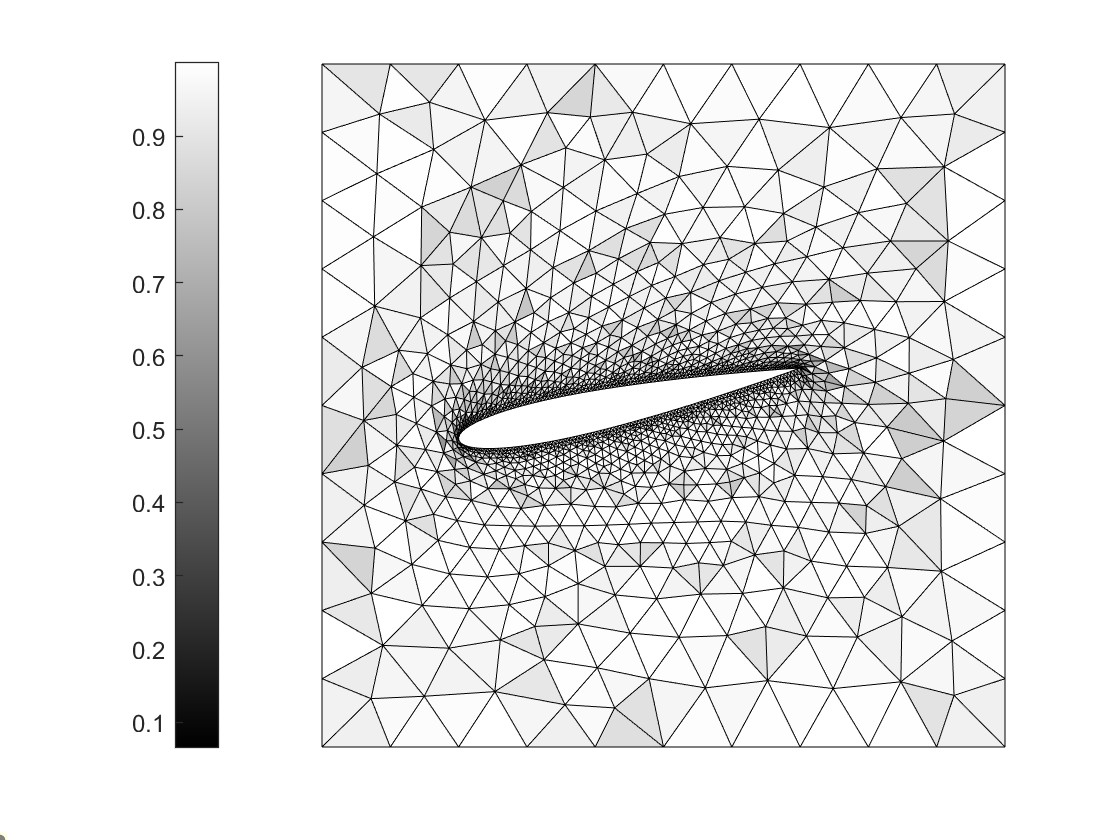}
        \caption{}
        \label{fig:NCAA0012MeshQualityTrue}
    \end{subfigure}
    \hfill
    \begin{subfigure}[t]{0.48\textwidth}
        \centering
        \includegraphics[width=\textwidth]{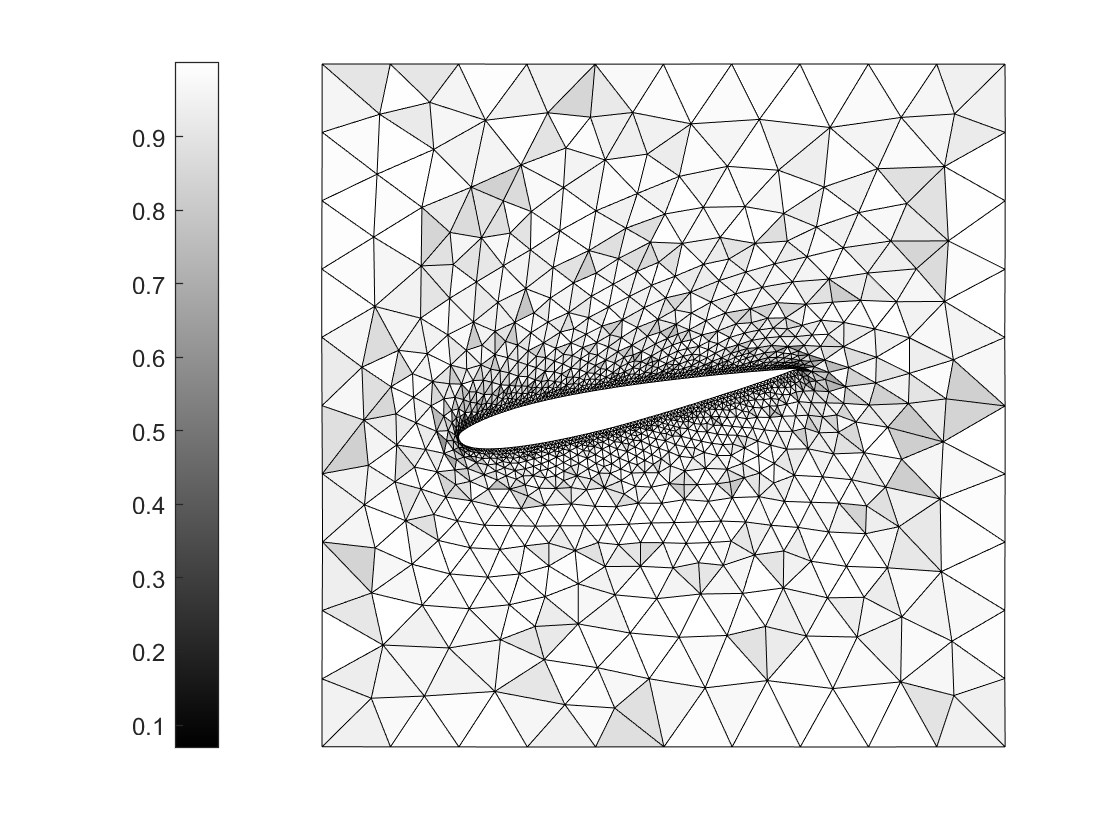}
        \caption{}
        \label{fig:NCAA0012MeshQualityPredict}
    \end{subfigure}
    \caption{Element quality metrics of NACA 0012 airfoil under combined translational and rotational rigid-body motion: (a) deformed mesh generated by FEM. (b) deformed mesh generated by the proposed surrogate model.}
\end{figure}

\begin{figure}[H]
    \centering
    \begin{subfigure}[t] {1\textwidth}
        \centering
        \includegraphics[width=1\textwidth]{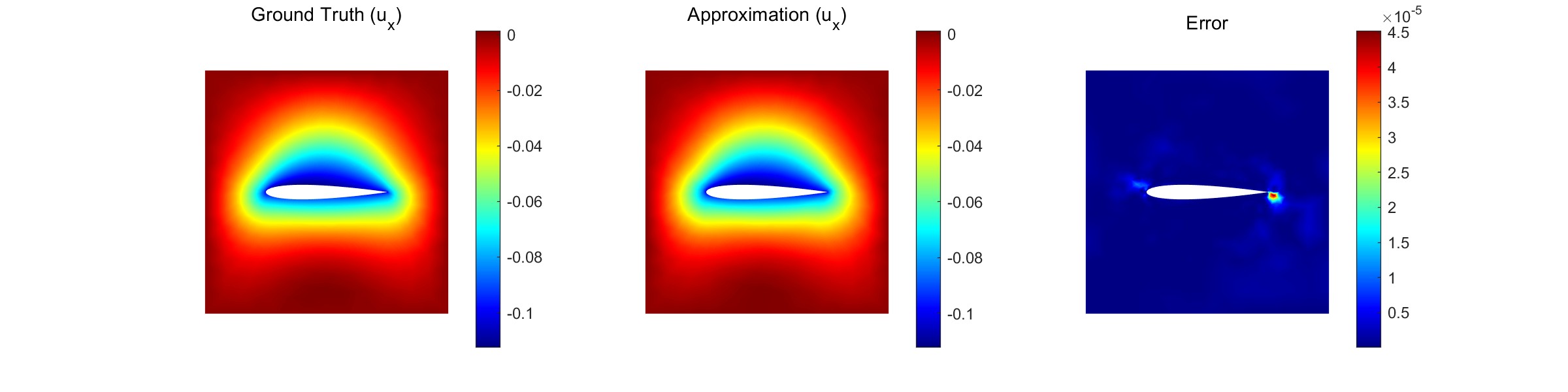} 
        \caption{}
        \label{fig:NCAA0012Ux}
    \end{subfigure}

    \begin{subfigure}[t] {1\textwidth}
        \centering
        \includegraphics[width=1\textwidth]{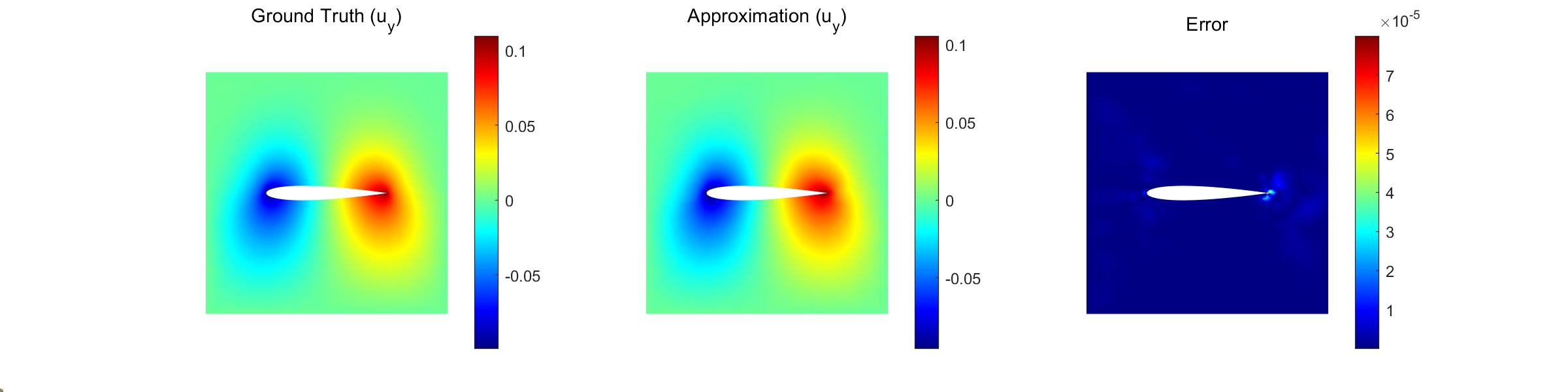} 
        \caption{}
        \label{fig:NCAA0012Uy}
    \end{subfigure}
    \caption{The proposed surrogate model for the solution of Eq. (\ref{eq:Linear_elastic}) in this case : (a) the solution in $X$ direction. (b) the solution in $Y$ direction. Both evaluate the approximation accuracy by the point-wise squared error. }
\end{figure}

The second case study examines a NACA 0012 airfoil profile defined by a nonlinear geometric boundary \cite{Huang2004}. The structure is centered within a computational domain of $(-0.5, 1.5) \times (-1, 1)$. The prescribed displacement field consists of a combined rigid body transformation, including a horizontal translation of $-0.1$ (leftward), a vertical shift of $+0.01$, and a counterclockwise rotation of $13^\circ$ about its geometric center. The resulting configuration, illustrating the combined effect of these boundary conditions, is shown in Fig. \ref{fig:NCAA0012DeformedMesh}.
The numerical results for the NACA 0012 airfoil demonstrate the robustness of the proposed method under combined translational and rotational transformations. As shown in Fig. \ref{fig:NCAA0012MeshQualityPredict}, the surrogate model maintains high element quality, which is visually indistinguishable from the FEM reference shown in Fig. \ref{fig:NCAA0012MeshQualityTrue}.
Notably, the \textbf{global relative error for this case is only 0.74\%}, a substantial reduction compared to the previous beam study. This improved performance is primarily attributed to the \textbf{geometric smoothness} of the airfoil profile. Unlike the rectangular beam, the NACA 0012 geometry lacks sharp corners, thereby avoiding the numerical singularities associated with boundary integral calculations at points of non-differentiability. Furthermore, the increased number of boundary sampling points required to capture the nonlinear curvature provides a richer set of constraints for the surrogate model, resulting in superior approximation accuracy in both the $U_x$ (Fig. \ref{fig:NCAA0012Ux}) and $U_y$ (Fig. \ref{fig:NCAA0012Uy}) displacement fields.
These findings suggest that the smoothness of the boundary and the density of its discretization are key factors in minimizing the approximation error of the proposed model.

\subsection{Verification of linear property}
\label{subsec:linearity}
\begin{figure}[H]
    \centering
    \begin{subfigure}[t]{0.48\textwidth}
        \centering
        \includegraphics[width=\textwidth]{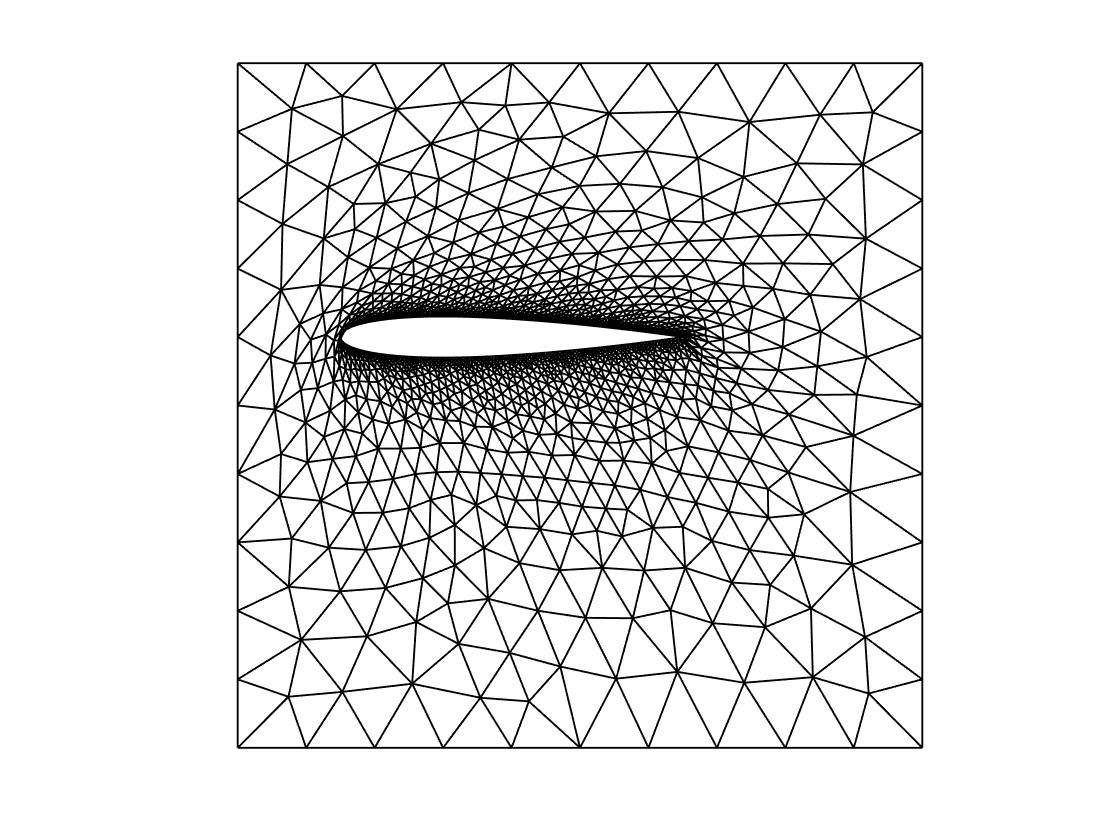}
        \caption{}
        \label{fig:scaleDeformedMesh}
    \end{subfigure}
    \hfill
    \begin{subfigure}[t]{0.48\textwidth}
        \centering
        \includegraphics[width=\textwidth]{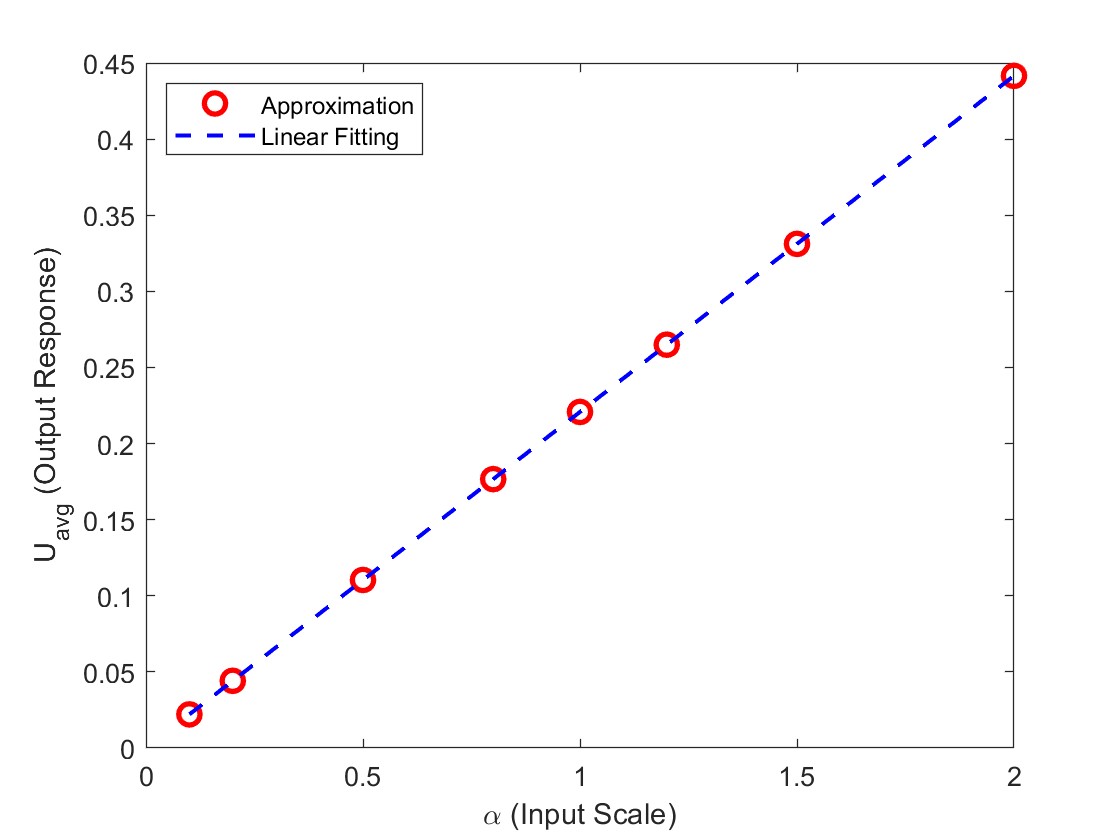}
        \caption{}
        \label{fig:scalePlot}
    \end{subfigure}
    \caption{The verification of the scale invariance of the NACA 0012 airfoil example: (a) The deformed mesh used as the scale reference. (b) The quantitative relationship between the scale scaling factor and output response. }
\end{figure}

\begin{figure}[H]
    \centering
    \begin{subfigure}[t]{0.48\textwidth}
        \centering
        \includegraphics[width=\textwidth]{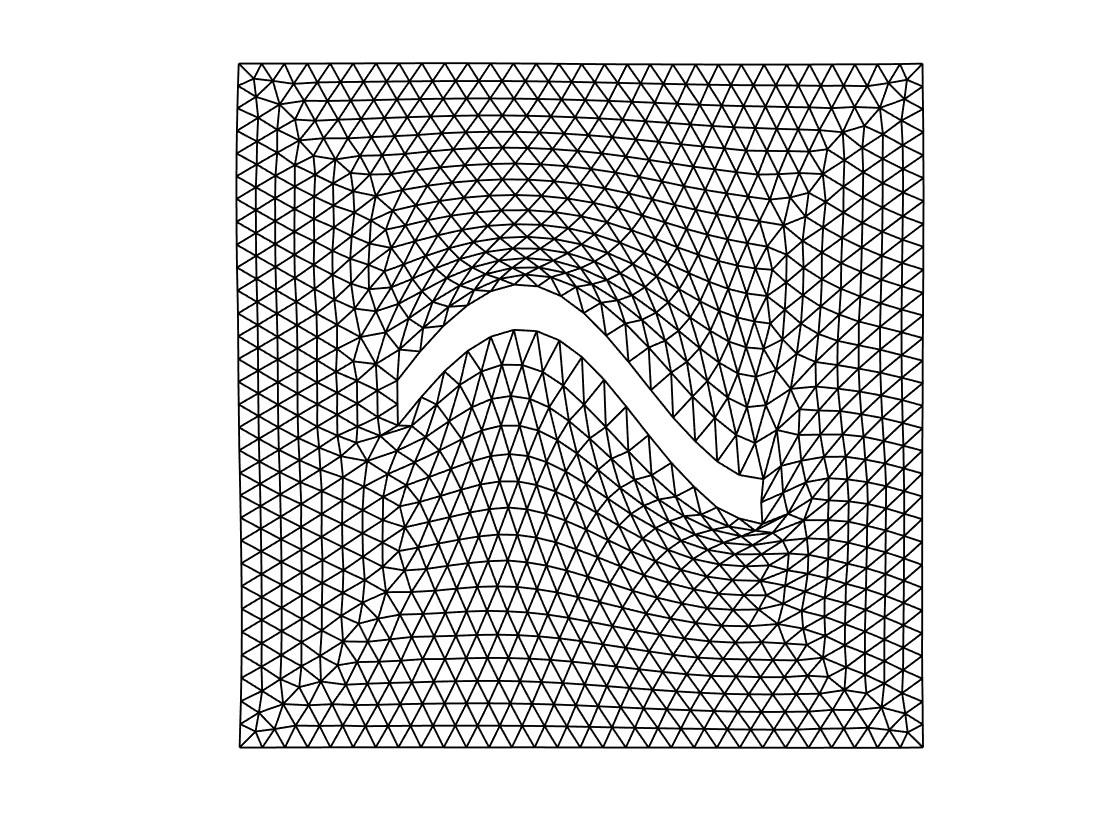}
        \caption{}
        \label{fig:scaleDeformedMesh1}
    \end{subfigure}
    \hfill
    \begin{subfigure}[t]{0.48\textwidth}
        \centering
        \includegraphics[width=\textwidth]{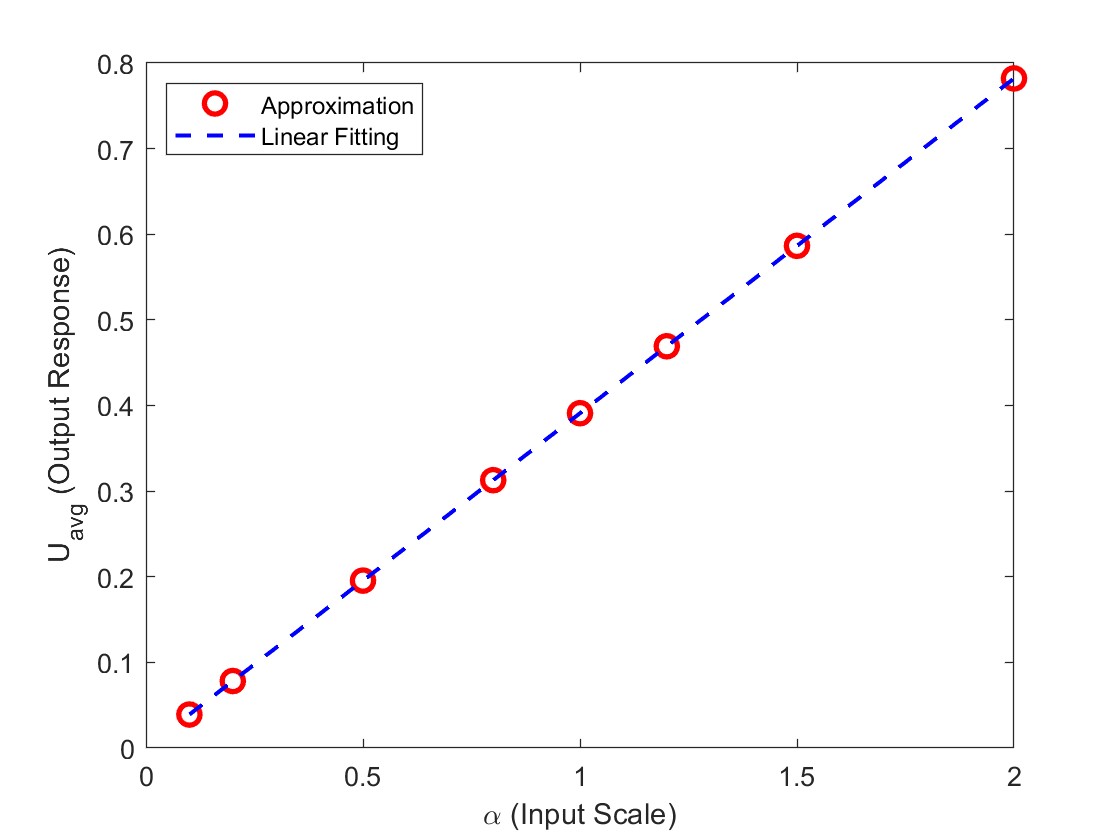}
        \caption{}
        \label{fig:scalePlot1}
    \end{subfigure}
    \caption{The verification of the scale invariance of the flexible beam example: (a) The deformed mesh used as the scale reference. (b) The quantitative relationship between the scale scaling factor and output response. }
\end{figure}

In this section, we verify the linearity of the proposed method. A mapping or operator is considered linear if it satisfies two fundamental principles: homogeneity and additivity (the superposition principle). To validate these properties in our surrogate model, we conduct tests on two distinct geometries: a NACA 0012 airfoil and a flexible beam. For the homogeneity test (scale invariance), the flexible beam is subjected to a base transverse bending mode defined by a sinusoidal displacement $\mathbf{d}_{base} = [0, \epsilon_y \sin(3\pi(x-3.5)/16)]^T$ with $\epsilon_y = 2.2$, while the NACA 0012 airfoil undergoes a base rigid-body transformation with translation parameters $t_x = -0.2$ and $t_y = 0.2$. For the superposition test, we focus on the NACA 0012 airfoil by defining two independent base modes: Mode A, representing a pure rotation of $\theta = -20^\circ$ around the geometric center $(0.5, 0)$, and Mode B, representing a pure translation with $t_x = 0.1$ and $t_y = 0.1$.The homogeneity property requires that scaling the input parameter by a factor $k$ results in a proportional scaling of the output response. We define the \textbf{Relative Linearization Error (RLE)} to quantify the deviation from this linear relationship:
\begin{equation}\text{RLE} = \frac{| \mathbf{U}p - k \cdot \mathbf{U}{base} |F}{k\cdot | \mathbf{U}{base} |F}
\end{equation}
where $\mathbf{U}_p$ is the displacement matrix generated by the surrogate model under the scaled input, $\mathbf{U}{base}$ is the reference displacement matrix at the standard scale ($k=1$), $k$ is the scaling factor, and $\| \cdot \|_F$ denotes the Frobenius norm. Similarly, the superposition principle states that the response caused by the sum of multiple inputs should equal the sum of the responses caused by each individual input. To quantify this, we define the \textbf{Relative Superposition Error (RSE)} as:
\begin{equation}
\text{RSE} = \frac{| \mathbf{U}_{comb} - (a \cdot \mathbf{U}_A + b \cdot \mathbf{U}_B) |F}{| \mathbf{U}{comb} |F}
\end{equation}
where $\mathbf{U}{comb}$ is the displacement matrix output by the model when multiple transformations are applied simultaneously; $\mathbf{U}_A$ and $\mathbf{U}_B$ are the displacement matrices for the independent base modes (Mode A and Mode B) described above; $a$ and $b$ are the combination coefficients representing the magnitude of each transformation; and $\| \cdot \|_F$ provides a global measure of error across all grid points in the mesh.
In the quantitative analysis of homogeneity, the output response (denoted as $U_{avg}$ in Figure \ref{fig:scalePlot} and \ref{fig:scalePlot1}) is defined as the mean Euclidean norm of the displacement vector across all $N$ grid points in the mesh. Specifically, for a given scaling factor $k$, the response is calculated as:
\begin{equation}
U_{avg} = \frac{1}{N} \sum_{i=1}^{N} \sqrt{u_i^2 + v_i^2}
\end{equation}
where $u_i$ and $v_i$ represent the predicted displacements in the $x$ and $y$ directions for the $i$-th node, respectively. 

As shown in Figure \ref{fig:scaleDeformedMesh} and Figure \ref{fig:scaleDeformedMesh1}, we applied various scaling factors $k$ to the reference deformed meshes of the NACA 0012 airfoil and the flexible beam, respectively. The quantitative results are presented in Figure \ref{fig:scalePlot} and Figure \ref{fig:scalePlot1}, where a strictly linear correlation between the scaling factor and the output response is observed.
\begin{table}[H]
\centering
\caption{Relative Linearization Error (RLE) for Different Scale Parameter Combinations in the NACA 0012 airfoil}
\label{tab:Scale_results}
\begin{tabular}{@{}lc @{\hspace{1cm}} lc@{}} 
\toprule
Scale $k$ & RLE & Scale $k$ & RLE \\ \midrule
0.1   & $3.1335 \times 10^{-7}$ & 1.2   & $3.3118 \times 10^{-7}$ \\
0.5   & $0.0000 \times 10^{0}$  & 1.5   & $3.2918 \times 10^{-6}$ \\
0.8   & $3.1335 \times 10^{-7}$ & 2.0   & $0.0000 \times 10^{0}$  \\
1.0   & $0.0000 \times 10^{0}$  & 2.5   & $3.3552 \times 10^{-7}$ \\ \bottomrule
\end{tabular}
\end{table}
\begin{table}[H]
\centering
\caption{Relative Linearization Error (RLE) for Different Scale Parameter Combinations in the flexible beam}
\label{tab:Scale_results2}
\begin{tabular}{@{}lc @{\hspace{1cm}} lc@{}} 
\toprule
Scale $k$ & RLE & Scale $k$ & RLE \\ \midrule
0.1   & $1.5924 \times 10^{-7}$ & 1.2   & $1.4415 \times 10^{-7}$ \\
0.5   & $0.0000 \times 10^{0}$  & 1.5   & $1.6150 \times 10^{-6}$ \\
0.8   & $1.5924 \times 10^{-7}$ & 2.0   & $0.0000 \times 10^{0}$  \\
1.0   & $0.0000 \times 10^{0}$  & 2.5   & $1.6491 \times 10^{-7}$ \\ \bottomrule
\end{tabular}
\end{table}
Table \ref{tab:Scale_results} and Table \ref{tab:Scale_results2} provide the calculated RLE for different combinations of $k$. In both cases, the RLE remains on the order of $10^{-7}$ or lower, which is effectively zero within machine precision. These results confirm that the model accurately preserves scale invariance across a wide range of parameters.

To further verify the linearity, we examine the superposition principle, which states that the response caused by the sum of multiple inputs should equal the sum of the responses caused by each input individually. For this test, we superimposed rotation and translation transformations on the NACA 0012 airfoil using parameters $a$ and $b$.
\begin{figure}[H]
    \centering
    \begin{subfigure}[t]{0.48\textwidth}
        \centering
        \includegraphics[width=\textwidth]{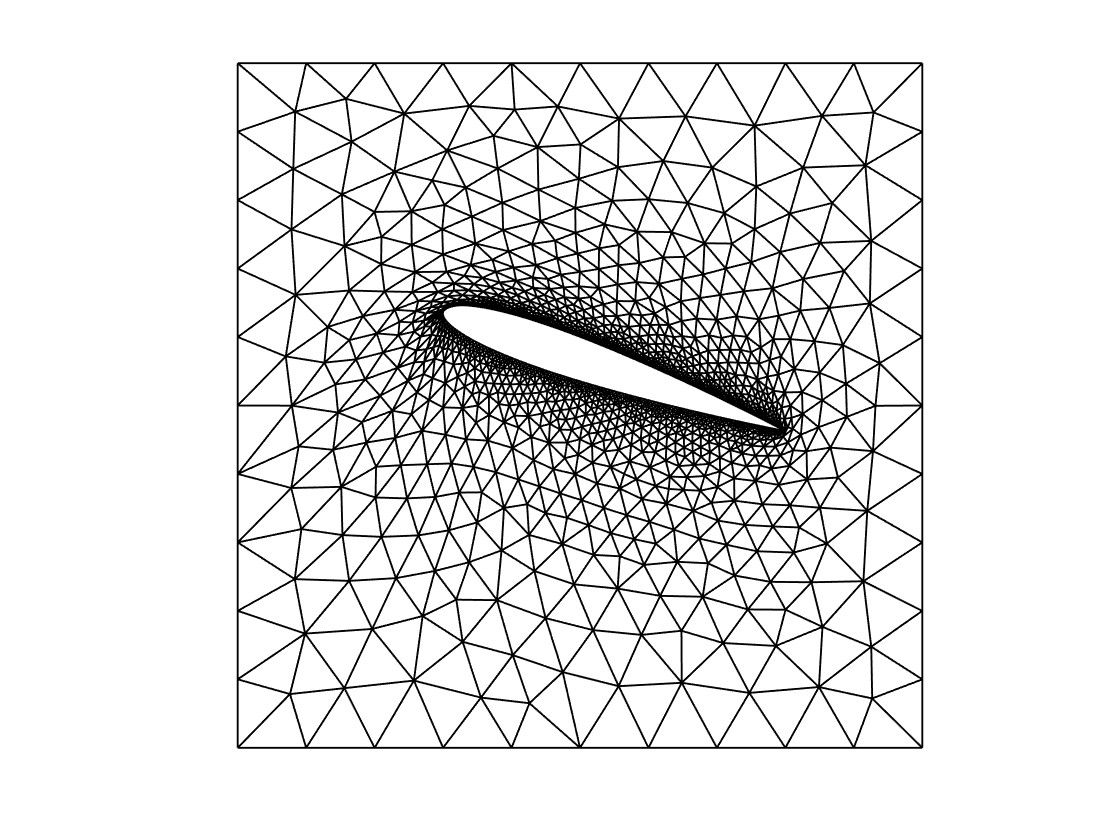}
        \caption{}
        \label{fig:addDeformedMesh}
    \end{subfigure}
    \hfill
    \begin{subfigure}[t]{0.48\textwidth}
        \centering
        \includegraphics[width=\textwidth]{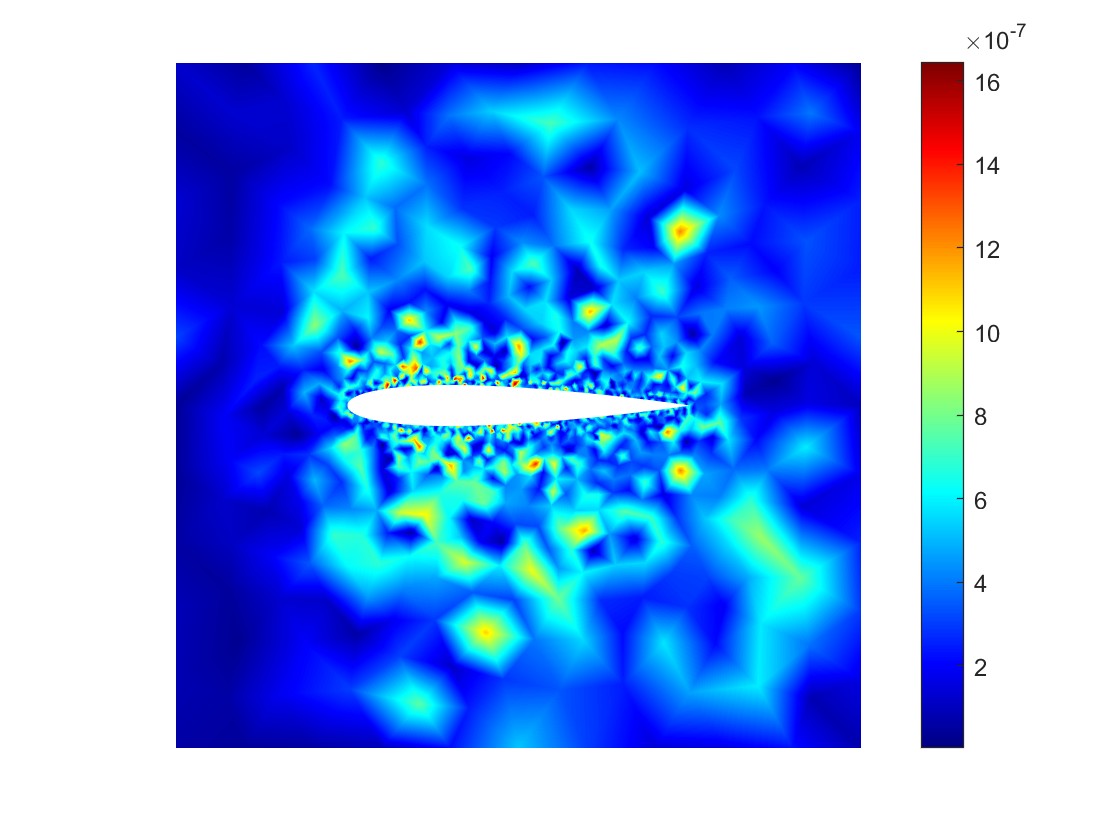}
        \caption{}
        \label{fig:addErrorDistribution}
    \end{subfigure}
    \caption{The verification of the superposition principle of the surrogate model: (a) is the deformed mesh generated by the NACA 0012 airfoil after rotation and translation superposition under the standard scale (a=1, b=1). (b) is the superposition error under the standard scale.}
\end{figure}
Figure \ref{fig:addDeformedMesh} illustrates the deformed mesh generated under a standard scale ($a=1, b=1$), while Figure \ref{fig:addErrorDistribution} displays the corresponding superposition error distribution. The spatial distribution of the error is negligible, indicating high fidelity in the combined transformation.
\begin{table}[H]
\centering
\caption{Relative Superposition Error (RSE) for Different Superposition Parameter Combinations}
\label{tab:add_rse_results}
\begin{tabular}{@{}lccc@{}}
\toprule
Combination & Parameter $a$ & Parameter $b$ & RSE \\ \midrule
1           & 1.0           & 1.0           & $3.7663 \times 10^{-6}$ \\
2           & 2.0           & 0.5           & $4.4212 \times 10^{-7}$ \\
3           & -1.0          & 1.0           & $3.0857 \times 10^{-7}$ \\
4           & 0.5           & -1.5          & $3.2926 \times 10^{-7}$ \\
5           & 0.1           & 2.0           & $3.6614 \times 10^{-7}$ \\
6           & 2.5           & 2.5           & $3.7587 \times 10^{-6}$ \\
7           & -1.2          & -1.2          & $3.7569 \times 10^{-6}$ \\
8           & 0.0           & 1.0           & $0.0000 \times  
10^{0\phantom{-}}$  \\ \bottomrule
\end{tabular}
\end{table}
The Relative Superposition Error (RSE) for various parameter combinations is summarized in Table \ref{tab:add_rse_results}. Across all tested combinations—including negative values and zero—the RSE consistently stays below $4 \times 10^{-6}$. The near-zero error across these diverse cases demonstrates that the surrogate model strictly adheres to the superposition principle.

The empirical evidence provided by the scaling tests and the superposition analysis confirms that the proposed surrogate model possesses robust linear properties. The extremely low values of RLE and RSE across different geometries and parameter ranges validate the mathematical consistency and reliability of the proposed method for complex mesh deformations.

\section{Conclusion and future Work}
In this study, we introduce a novel mesh deformation framework that integrates boundary integral theory with neural operators. Employing deep neural networks to learn the Green’s traction kernel of linear elasticity equations, we effectively address the challenges of mesh deformation within a boundary value problem (BVP) context. Numerical experiments on flexible beams and NACA0012 airfoils, along with rigorous linearity verification, demonstrate that our model achieves robust generalization across diverse boundary conditions while strictly adhering to physical principles. In the future, several research directions will be pursued to improve the utility of the framework. First, we aim to extend the model from 2D to 3D deformation by incorporating high-efficiency computational techniques to manage the increased complexity. Second, we plan to integrate specific geometric descriptors to enable cross-geometry generalization across diverse structural topologies. Furthermore, we will refine the architecture to capture nonlinear elastic behaviors, such as hyperelasticity, enabling rapid simulations of large deformations. Finally, we intend to incorporate physics-informed constraints to reduce the model's reliance on large-scale training data, thereby improving its efficiency and robustness in data-scarce engineering environments.

\bibliographystyle{unsrt}  
\bibliography{references}

\end{document}